\documentclass{ccgrt} 

\newcommand{\ideal}{\trianglelefteq}
\newcommand{\coideal}{\dot{\trianglelefteq}}
\newcommand{\coleq}{\dot{\leq}}

\newcommand{\llongrightarrow}
{\relbar\joinrel\relbar\joinrel\relbar\joinrel\rightarrow}
\newcommand{\loongtwoheadrightarrow}
{\relbar\joinrel\relbar\joinrel\twoheadrightarrow}

\newcommand{\gerg}{\mathfrak{g}}

\newcommand{\gerk}{\mathfrak{k}}

\newcommand{\gert}{\mathfrak{t}}

\newcommand{\gerI}{\mathfrak{I}}
\newcommand{\gerC}{\mathfrak{C}}
\newcommand{\calI}{\mathcal{I}}
\newcommand{\calC}{\mathcal{C}}
\newcommand{\calS}{\mathcal{S}}

\newcommand{\h}{\hbar}
\newcommand{\HA}{{\mathcal{HA}}}
\newcommand{\Hpicc}{{\scriptscriptstyle \mathcal{H}}}
\newcommand{\Ppicc}{{\scriptscriptstyle \mathcal{P}}}

\newcommand{\N}{\mathbb{N}}

\newcommand{\kh}{\Bbbk [[\h]]}
\newcommand{\id}{\text{\sl id}}

\newcommand{\otimeshat}{\,\widehat{\otimes}\,}
\newcommand{\otimestilde}{\,\widetilde{\otimes}\,}

\newcommand{\QUEA}{\mathcal{Q\hskip1ptUE\hskip-1ptA}}

\newcommand{\QFSHA}{\mathcal{Q\hskip1ptFSHA}}

\newcommand{\uhg}{U_\h(\gerg)}
\newcommand{\fhg}{F_\h[[G]]}

\begin{document}

\title{Quantum duality principle for coisotropic  \\
  subgroups and Poisson quotients\,\footnote{\,{\it
\uppercase{M}\uppercase{S}\uppercase{C}\,2000:\,}
Primary 17\uppercase{B}37, 20\uppercase{G}42, 58\uppercase{B}32;
Secondary 81\uppercase{R}50.
\newline {\it \uppercase{K}eywords:\,}\
quantum groups, \uppercase{P}oisson homogeneous
spaces, coisotropic subgroups.}}

\author{Nicola Ciccoli}
\address{Universit\`a di Perugia,
Dipartimento di Matematica,  \\
  Via Vanvitelli 1, I-06123 Perugia --- ITALY
 \vspace{1mm}  \\
E-mail: ciccoli@dipmat.unipg.it }

\author{Fabio Gavarini}
\address{Universit\`a di Roma ``Tor Vergata'',
Dipartimento di Matematica,  \\
  Via della ricerca scientifica 1, I-00133 Roma  ---  ITALY
 \vspace{1mm}  \\
E-mail: gavarini@mat.uniroma2.it }

\maketitle

\abstracts
  {We develop a quantum duality principle for coisotropic subgroups of
a (formal) Poisson group and its dual: namely, starting from a quantum
coisotropic subgroup (for a quantization of a given Poisson group)
we provide functorial recipes to produce quantizations of the dual
coisotropic subgroup (in the dual formal Poisson group).  By the natural
link between subgroups and homogeneous spaces, we argue a quantum
duality principle for Poisson homogeneous spaces which are Poisson
quotients, i.e.~have at least one zero-dimensional symplectic leaf.
                                             \\
   Only bare results are presented, while detailed proofs can be
found in [3].}

\setcounter{section}{0}

\section{Introduction}

\vskip5pt

\noindent
 In the study of quantum groups, the natural semiclassical counterpart
is the theory of deformation (or quantization) of Poisson groups: actually,
Drinfeld himself introduced Poisson groups as the semiclassical limits of
quantum groups.  Therefore, it should be not surprising that the geometry
of quantum groups turns more clear and comprehensible when its connection
with Poisson geometry is more transparent.  The same situation occurs
when dealing with Poisson homogeneous spaces of Poisson groups.
                                             \par \noindent
 In particular, in the study of Poisson homogeneous spaces a special
role is played by  {\sl Poisson quotients}.  By this we mean
Poisson homogeneous spaces whose symplectic foliation has at least
one zero-dimensional leaf: therefore, they can be seen as pointed 
Poisson homogeneous spaces, just like Poisson groups themselves
are pointed by the identity element.
                                    \par \noindent
 Poisson quotients form a natural subclass of Poisson homogeneous
$ G $--spaces  ($ G $  a Poisson group) which is best adapted to
the standard relation between homogeneous  $ G $--spaces  and
subgroups of  $ G \, $:  to a given Poisson quotient, one
associates the stabilizer subgroup of its distinguished point (the
 fixed zero-dimensional symplectic leaf).  What characterizes such subgroups is
{\sl coisotropy},  with respect to the Poisson structure on  $ G $
(see the definition in Section 3 later on).  On the other hand, if
a (closed) subgroup  $ K $  of  $ G $  is coisotropic, then the
homogeneous $ G $--space  $ G \big/ K $  is a Poisson quotient.
So the two notions of Poisson quotient and coisotropic subgroup
must be handled in couple.  In particular, the quantization
process for a Poisson $ G $--quotient corresponds to a similar
procedure for    
   \hbox{the attached coisotropic subgroup of  $ G $.}
                                    \par \noindent
 If one looks at quantizations of a Poisson homogeneous
space, their existence is guaranteed only if the space is
a  quotient [8];  thus the notion of Poisson quotient
shows up naturally as a necessary condition.
On the other hand, let  $ K $  be a subgroup of  $ G \, $, 
and assume that  $ G $  has a quantization, inducing on
it a Poisson group structure.  If  $ K $  itself also admits a
quantization, which is ``consistent'' (in a natural sense) with
the one of  $ G \, $,  then  $ K $  is
automatically coisotropic in  $ G \, $.  So also the related
notion of coisotropic subgroup  shows to be
a necessary condition for the existence of quantizations.
Of course an analogous description can be entirely carried out at an infinitesimal level, with conditions
at the level of Lie bialgebras.
                                    \par \noindent
 When dealing with quantizations of Poisson groups (or Lie bialgebras), a
precious tool is the quantum duality principle (QDP).  Roughly speaking,
it claims that any quantized enveloping algebra can be turned   --- via
a functorial recipe ---   into a quantum function algebra for the dual
Poisson group; conversely, any quantum function algebra can be turned
into a quantization of the enveloping algebra of the dual Lie bialgebra. 
To be precise, let  $ \QUEA $  and  $ \QFSHA $ respectively be the
category of all quantized universal enveloping algebras (QUEA) and
the category of all quantized formal series Hopf algebras (QFSHA),
in Drinfeld's sense.  Then the QDP establishes [6,11] a category
equivalence between  $ \QUEA $  and  $ \QFSHA $  via two functors, 
$ \, (\ )' \colon \QUEA \longrightarrow \QFSHA \, $  and  $ \,
(\ )^\vee \colon \QFSHA \longrightarrow \QUEA \, $.  Moreover,
starting from a QUEA over a Lie bialgebra (resp.~from a QFSHA over
a Poisson group) the functor  $ (\ )' $  (resp.~$ (\ )^\vee \, $) 
gives a QFSHA (resp.~a QUEA) over the dual Poisson group (resp.~the
dual Lie bialgebra).  In short,  $ {U_\h(\gerg)}' \!\! = F_\h[[G^*]]
\, $  and  $ {F_\h[[G]]}^\vee \!\! = U_\h(\gerg^*) \, $  for any Lie
bialgebra  $ \gerg $  and Poisson group  $ G $  with  $ \, \text{\it
Lie}(G) = \gerg \, $.  So from a quantization of any Poisson group
this principle gets out a quantization of the dual Poisson group too.  
                                                 \par \noindent
 In this paper we establish a similar quantum duality principle for
(closed) coisotropic subgroups of a Poisson group  $ G $,  or equivalently for Poisson  $ G $--quotients,  sticking to the formal approach (hence dealing with quantum groups \`a la Drinfeld).  The starting point is that any formal coisotropic subgroup  $ K $  of a Poisson group  $ G $  has two possible algebraic descriptions via objects related to  $ U(\gerg) $  or  $ F[[G]] $,  and similarly for the formal Poisson quotient  $ G\big/K \, $;  thus the datum of  $ K $  or equivalently of  $ G\big/K $  is described algebraically in four
possible ways.  By  {\sl quantization\/}  of such a datum we mean
a quantization of any one of these four objects, which has to be ``consistent''   --- in a natural sense ---   with given quantizations  $ U_\h(\gerg) $  and  $ \fhg $  of  $ G \, $.  Our ``QDP'' now is a bunch of functorial recipes to produce, out of a quantization of
$ K $  or  $ G\big/K $  as before, a similar quantization of the
so-called  {\sl complementary dual\/}  of  $ K \, $,  that is the coisotropic subgroup  $ K^\perp $  of  $ G^* $  whose tangent Lie bialgebra is just  $ \gerk^\perp $  inside  $ \gerg^* \, $,  or of
the associated Poisson  $ G^* $--quotient,  namely  $ G^* \big/
K^\perp \, $.  The basic idea is quite simple.  The quantizations of coisotropic subgroups   --- or Poisson quotients ---   are sub-objects of quantizations of Poisson groups, and the recipes of the original QDP (for Poisson groups) apply to the latter objects.  Then we simply ``restrict'', somehow, such recipes to the previously mentioned sub-objects.
                                                 \par \noindent
 In recent times, the general problem of quantizing coisotropic manifolds
 of a given Poisson manifold, in the context of deformation quantization,
has raised quite some interest [1,3].  It is then important to point
out that ours is by no means an existence result: instead, it can be
thought of as a  {\sl duplication result},  because it yields a new
quantization   --- for a complementary dual object ---   out of one
given from scratch (much like the QDP for quantum groups).  On the
other hand, we would better stress that our result is really  {\sl
effective},  and calling for applications.  A sample of application
is presented in the extended version of this work [3]; see also
Subsection 5.6.

\section{The classical setting}

\vskip5pt

\subsection{Formal Poisson groups}    
 We shall work in the setup of
formal geometry.  Recall that a formal variety is uniquely characterized
by a tangent or cotangent space (at its unique point), and it is described
by its ``algebra of regular functions''   --- such as  $ F[[G]] $  below.
This is a complete, topological local ring which can be realized
as a  $ \Bbbk $--algebra  of formal power series.  Hereafter  $ \Bbbk $
is a field of zero characteristic.
                                            \par \noindent
 Let  $ \gerg $  be a finite dimensional Lie algebra over  $ \Bbbk
\, $,  and let  $ U(\gerg) $  be its universal enveloping algebra (with
the natural Hopf algebra structure).  We denote by  $ F[[G]] $  the
algebra of functions on the formal algebraic group  $ G $  associated
to  $ \gerg $  (which depends only on  $ \gerg $  itself); this is
a complete, topological Hopf algebra.  Furthermore  $ \, F[[G]] \cong
{U(\gerg)}^* \, $,  so that there is a natural pairing of (topological)
Hopf algebras   --- see below ---   between  $ U(\gerg ) $  and
$ F[[G]] \, $.
                                            \par \noindent
 In general, if  $ H $,  $ K $  are Hopf algebras (even topological)
over a ring  $ R \, $,  a pairing  $ \, \langle \,\ , \,\ \rangle \,
\colon \, H \times K \longrightarrow R \, $   is called a  {\sl Hopf
pairing\/}  if
  $ \; \big\langle x, y_1 \cdot y_2 \big\rangle = \big\langle \Delta(x),
y_1 \otimes y_2 \big\rangle \, ,  \;\;  \big\langle x_1 \cdot x_2,
y \big\rangle = \big\langle x_1 \otimes x_2, \Delta(y) \big\rangle
\, $,  $ \; \langle x, 1 \rangle = \epsilon(x) \, $,  $ \; \langle 1,
y \rangle = \epsilon(y) \, $, $ \; \big\langle S(x), y \big\rangle =
\big\langle x, S(y) \big\rangle \; $  for all  $ \, x , x_1 , x_2 \in
H \, $,  $ \, y, y_1, y_2 \in K \, $.  The pairing is called
{\sl perfect\/}  if it is non-degenerate.
                                            \par \noindent
Assume  $ G $  is a formal  {\sl Poisson\/}  (algebraic) group.
Then  $ \gerg $  is a Lie bialgebra,  $ U(\gerg) $  is a co-Poisson Hopf
algebra,  $ F[[G]] $  is a topological Poisson Hopf algebra, and
the Hopf pairing above respects these additional co-Poisson and Poisson
structures.  Furthermore, the linear dual  $ \gerg^* $  of  $ \gerg $
is a Lie bialgebra as well, so a dual formal Poisson group  $ G^* $
exists.

\vskip2pt

 $ \underline{\hbox{\sl Notation}} $:  hereafter, the symbol  $ \,
\coideal \, $  stands for ``coideal'',  $ \, \leq^1 \, $  for ``unital
subalgebra'',  $ \, \coleq \, $  for ``subcoalgebra'',  $ \, \leq_\Ppicc
\, $  for ``Poisson subalgebra'',  $ \, \coideal_\Ppicc \, $  for
``Poisson coideal'',  $ \, \leq_\Hpicc \, $  for ``Hopf subalgebra'',
$ \, \ideal_\Hpicc \, $  for ``Hopf ideal'', and the subscript  $ \ell $
stands for ``left'' (everything in topological sense if necessary).

\subsection{Subgroups and homogeneous  $ G $--spaces} 
 A homogeneous left  $ G $--space  $ M $  corresponds to a conjugacy
class of closed subgroups  $ \, K = K_M \, $,  \, which we assume
connected, of  $ G \, $,  such that  $ \, M \cong G \big/ K \, $. 
In formal geometry  $ K $  may be replaced by  $ \, \gerk :=
\text{\it Lie}(K) \, $.  The whole geometric setting given by
the pair  $ \big( K, G/K \big) $  then is encoded by any one
of the following data:

\vskip3pt

   {\it (a)} \  the set  $ \, \calI = \calI(K) \equiv \calI(\gerk) \, $
of (formal) functions vanishing on  $ K \, $,  that is to say  $ \,
\calI = \big\{ \varphi \!\in\! F[[G]] \,\big|\, \varphi(K) \! = \! 0
\big\} \, $;  \, note that  $ \; \calI \ideal_\Hpicc \! F[[G]] \; $;
                                                 \par
   {\it (b)} \  the set of left  $ \gerk $--invariant  functions,
namely  $ \, \calC = \calC(K) \equiv \calC(\gerk) = {F[[G]]}^K \; $;
\, note that  $ \; \calC \leq^1 \! \coideal_\ell \, F[[G]] \; $;
                                                 \par
   {\it (c)} \  the set  $ \, \gerI = \gerI(K) \equiv \gerI(\gerk) \, $
of left-invariant differential operators on  $ F[[G]] $   which
vanish on  $ \, {F[[G]]}^K \, $,  \, that is  $ \, \gerI = U(\gerg)
\cdot \gerk \, $  (via standard identifi\-cation of the set of
left-invariant differential operators with  $ U(\gerg) \, $);
\, note that  $ \; \gerI(\gerk) = \gerI \ideal_\ell \!
\coideal \, U(\gerg) \, $;
                                                 \par
   {\it (d)} \  the universal enveloping algebra of  $ \, \gerk \, $,
\, denoted  $ \gerC = \gerC(K) \equiv \gerC(\gerk) :=  U(\gerk) \; $;
\, note that  $ \; \gerC(\gerk) = \gerC \leq_\Hpicc \! U(\gerg) \; $.

\vskip4pt

\noindent
 In this way any formal subgroup  $ K $  of  $ G \, $,  or the
associated homogeneous  $ G $--space  $ \, G \big/ K \, $,  \, is
characterized by any of the following objects:
  $$  \text{\it (a)} \hskip4pt \calI \ideal_\Hpicc \! F[[G]]
\hskip7pt  \text{\it (b)} \hskip4pt \calC \leq^1 \! \coideal_\ell
\, F[[G]]  \hskip7pt  \text{\it (c)} \hskip4pt \gerI \ideal_\ell
\coideal \, U(\gerg)  \hskip7pt  \text{\it (d)} \hskip4pt
\gerC \leq_\Hpicc \! U(\gerg)  $$
These four data are all equivalent to each other, as we now explain.
                                                 \par \noindent
 For any Hopf algebra  $ H \, $,  with counit  $ \epsilon \, $,  and
every submodule  $ \, M \subseteq H \, $,  \, we set:  $ \, M^+ := M
\cap \text{\it Ker}\,(\epsilon) \, $  and  $ \, H^{\text{\it co}M} :=
\big\{\, y \in H \,\big|\, \big( \Delta(y) - y \otimes 1 \big) \in H
\otimes M \,\big\} \, $  (the set of  $ M $--{\sl coinvariants\/}  of
$ H \, $).  Letting  $ \mathbb{A} $  be the set of all subalgebras left
coideals of  $ H $  and  $ \mathbb{K} $  be the set of all coideals
left ideals of  $ H \, $,  we have well-defined maps  $ \, \mathbb{A}
\longrightarrow \mathbb{K} \, $,  $ \, A \mapsto H \cdot A^+ \, $,
and  $ \, \mathbb{K} \longrightarrow \mathbb{A} \, $,  $ \, K \mapsto
H^{\text{\it co}K} \, $  (see for instance Masuoka's work [15]).  Then
the above equivalence stems from

 \vskip2pt

   {\it --- (1)}  \,  {\sl orthogonality relations\/}   --- w.r.t.~the
natural pairing between  $ F[[G]] $  and  $ U(\gerg) $  ---   namely
$ \, \calI = \gerC^\perp \, $,  $ \, \gerC = \calI^\perp \, $,  \;
and  $ \; \calC = \gerI^\perp \, $,  $ \, \gerI = \calC^\perp \, $;
                                                 \par
   {\it --- (2)}  \,  {\sl subgroup-space correspondence},
namely  $ \, \calI = F[[G]] \cdot \calC^+ \, $,  $ \, \calC
= {F[[G]]}^{\text{\it co}\calI} \, $,  \, and  $ \; \gerI =
U(\gerg) \, \gerC^+ \, $,  $ \, \gerC = {U(\gerg)}^{\text{\it
co} \gerI} \, $.  Moreover, the maps  $ \, \mathbb{A} \longrightarrow
\mathbb{K} \, $  and  $ \, \mathbb{K} \longrightarrow \mathbb{A} \, $
above are inverse to each other in the formal setting.

\subsection{Coisotropic subgroups and Poisson quotients}  
 Assume now that  $ G $  is a formal Poisson group.  A closed formal
subgroup  $ K $  of  $ G $  with Lie algebra  $ \gerk $  is called
{\sl coisotropic}  if its defining ideal  $ \calI(\gerk) $  is a
topological Poisson subalgebra of  $ F[[G]] \, $.  The following
are equivalent [13,14]:

\vskip1pt

   {\it (C-i)}  \,  $ K $  is a coisotropic formal subgroup of  $ G \, $;
                                                 \par
   {\it (C-ii)}  \,  $ \delta(\gerk) \subseteq \gerk \wedge \gerg \, $,
that is  $ \gerk $  is a Lie coideal of  $ \gerg \, $;
                                                 \par
   {\it (C-iii)}  \,  $ \gerk^\perp \, $  is a Lie subalgebra
of  $ \gerg^* \, $

\vskip3pt

\noindent
Clearly, conditions  {\it (C-i,ii,iii)\/}  characterize coisotropic subgroups.
                                                 \par \noindent
   As to homogeneous spaces, a formal Poisson manifold  $ \, (M,
\omega_M) \, $  is a  {\sl Pois\-son homogeneous  $ G $--space\/}
if it carries a homogeneous action  $ \, \phi \colon \, G \times
M \rightarrow M \, $   which is a smooth Poisson map.  
%
%
%
In addition,  $ \, (M,\omega_M) \, $  is said to be  {\sl of group
type\/}  (after Drinfeld [7]), or simply a  {\sl Poisson
quotient},  if there is a coisotropic closed Lie subgroup  $ K_M $ 
of  $ \, G $  such that  $ \, G \big/ K_M \simeq M \, $  and the
natural projection map  $ \, \pi \colon \, G \longrightarrow
G \big/ K_M \simeq M \, $  is a Poisson map.
                                                 \par \noindent
   The following is a characterization of Poisson quotients [17]:

\vskip3pt

   {\it (PQ-i)}  \,\; there exists  $ \, x_0 \in M \, $  whose
stabilizer  $ \, G_{x_0} \, $  is coisotropic in  $ \, G \, $;
                                                 \par
   {\it (PQ-ii)}  \,\, there exists  $ \, x_0 \in M \, $  such that
$ \, \phi_{x_0} \, \colon \, G \longrightarrow M \, $,  $ \, g \mapsto
\phi(g,x_0) \, $,  \; is a Poisson map, that is  $ M $  is a Poisson
quotient;
                                                 \par
   {\it (PQ-iii)}  \, there exists  $ \, x_0 \in M \, $  such that
$ \, \omega_M(x_0) = 0 \, $.

\vskip4pt

\noindent
It is important to remark that in Poisson geometry, the usual relationship between closed subgroups of  $ G $  and  $ G $--homogeneous spaces does not hold anymore: in fact, in the  {\sl same\/}  conjugacy class one can have Poisson subgroups, coisotropic subgroups  {\sl and\/}  even non-coisotropic subgroups.  Now, the above characterization means exactly that the Poisson quotients are just those Poisson homogeneous spaces in which (at least) one of the stabilizers is coisotropic. 
Moreover, in general the correspondence between homogeneous spaces and subgroups is somewhat ambiguous, because it passes through the choice of a distinguished point of the space (whose stabilizer is the subgroup).  In the case of Poisson quotients this ambiguity is cleared off, as we do fix as distinguished point on the space the zero-dimensional symplectic leaf that it has for sure   --- although it is non-unique, a priori.
                                            \par \noindent
 In addition, passing through coisotropic subgroups allows us to introduce a good notion of (Poisson)  {\sl duality\/}  for our objects, namely

\begin{definition} \label{coiso}
  \, (notation of Subsection 2.1)
                                          \par
   {\it (a)}  If  $ K $  is a formal coisotropic subgroup of
$ \, G $,  we call  {\sl complementary dual\/}  of  $ \, K $  the
formal subgroup  $ K^\perp $  of  $ \, G^* $  whose tangent Lie
algebra is  $ \, \gerk^\perp \, $.
                                          \par
   {\it (b)}  If  $ \, M \cong G \big/ K_M \, $  is a formal
Poisson  $ G $--quotient,  with  $ K_M $  coiso-tropic, we call
$ \, M^\perp := G^* \big/ K_{\!M}^{\;\perp} \, $  the  {\sl
complementary dual\/}  of  $ M $.
\end{definition}

\noindent    
 Here the key point is that   --- by  {\it (C-iii)\/}  in Subsection 2.3 ---   a subset  $ \gerk $  of  $ \gerg $  is a Lie coideal if and only if  $ \gerk^\perp $  is a Lie subalgebra of  $ \gerg^* \, $.  Even more, by  {\it (C-i,ii,iii)},  the complementary dual subgroup to a coisotropic subgroup is  {\sl coisotropic\/}  too, and taking twice the complementary dual gives back the initial subgroup.  Similarly, the Poisson homogeneous space which is complementary dual to a Poisson homogeneous space of group type is in turn  {\sl of group type\/}  too, and taking twice the complementary dual gives back the initial manifold.  At the level of Poisson homogeneous spaces, one should think of  $ \big( K, G/K^\perp \big) $  and  $ \big( G/K, K^\perp ) $  as mutually dual pairs; if  $ \, K = \{ e \} \, $,  one recovers the usual couple of Poisson groups $ G $,  $ G^* \, $.  When using these pairs a price is paid: one object (the subgroup) is not a Poisson manifold and the other (the homogeneous spaces) is not a group anymore.

\subsection{Remarks}
 {\it (a)}  The notion of Poisson homogeneous
$ G $--spaces of group type was introduced by 
Drinfeld [7], who also explained the relation
between such  $ G $--spaces  and Lagrangian subalgebras of
Drinfeld's double  $ \, D(\gerg) = \gerg \oplus \gerg^* \, $. 
This was further developed by Evens and Lu [10], who gave a
Poisson structure on the algebraic variety of Lagrangian
subalgebra. It is an open problem to quantize this sort
of  {\it universal moduli space\/}  of Poisson homogeneous 
$ G $--spaces.
                                            \par
   {\it (b)}  As a matter of notation, we denote by  $ \, \hbox{\it co}
\hskip1pt\calS(G) \, $  the set of all formal coisotropic subgroups of
$ G \, $,  which is as well described by the set of all Lie subalgebras,
Lie coideals of  $ \gerg \, $.  Since it is ordered by inclusion, this
set will be also considered as a category.   $ \quad \diamondsuit $

\subsection{Algebraic characterization of coisotropy}    
 Let  $ K $ be a formal co\-isotropic subgroup of  $ G $.  In terms
of Subsection 2.2, coisotropy corresponds to
 \vskip-16pt
  $$  \text{\it (a)} \;\; \calI \leq_\Ppicc F[[G]]  \hskip17pt
\text{\it (b)} \;\; \calC \leq_\Ppicc F[[G]]  \hskip17pt
\text{\it (c)} \;\; \gerI \, \coideal_\Ppicc \, U(\gerg)  \hskip17pt
\text{\it (d)} \;\, \gerC \, \coideal_\Ppicc \, U(\gerg)  $$
 \vskip-4pt
\noindent
 Thus a formal coisotropic subgroup of  $ G $  is  identified by any
one of
 \vskip-13pt
  $$  \begin{matrix}
   \text{\it (a)} \;\;\;\; \calI \ideal_\Hpicc \leq_\Ppicc \! F[[G]]
&  \hskip23pt  \phantom{\Big|}  \text{\it (b)} \;\;\;\; \calC \leq^1
\! \coideal_\ell \leq_\Ppicc \! F[[G]]  \\
   \text{\it (c)} \;\;\;\; \gerI \ideal_\ell \! \coideal \,
\coideal_\Ppicc \, U(\gerg)  &  \hskip23pt  \phantom{\Big|}
\text{\it (d)} \;\;\;\; \gerC \leq_\Hpicc \! \coideal_\Ppicc
\, U(\gerg)
\end{matrix}  $$

\vskip3pt

\section{The quantum setting}

\vskip5pt

\subsection{Topological  $ \Bbbk[[\h]] $--modules  and tensor structures} 
 Let  $ \kh $  be the topological ring of formal power series in the
indeterminate  $ \h \, $.  If  $ X $  is any  $ \kh $--module, we set
$ \, X_0 := X \big/ \h X = \Bbbk \otimes_{\kh} X \, $, \, the  {\sl
specialization\/}  of  $ X $  at  $ \, \h = 0 \, $,  or  {\sl
semiclassical limit\/}  of  $ X \, $.  We are interested in certain
families of  $ \kh $--modules,  complete with respect to suitable
topologies, and for which a good notion of tensor product be available.
The exact notion of topology or of tensor product to choose depends on
whether one looks for quantizations of universal enveloping algebras or
of formal series Hopf algebras: this leads to two different setups, as
follows.
                                            \par \noindent
 First, let  $ \mathcal{T}_{\otimeshat} $  be the category whose objects
are all topological  $ \kh $--modu\-les  which are topologically free and
whose morphisms are the  $ \kh $--linear  maps (which are automatically
continuous).  It is a tensor category for the tensor product  $ \,
T_1 \otimeshat T_2 \, $  defined as the separated  $ \h $--adic
completion of the algebraic tensor product  $ \, T_1 \otimes_{\kh}
T_2 \, $  (for all  $ T_1 $,  $ T_2 \in \mathcal{T}_{\otimeshat} $).
We denote by  $ \HA_{\otimeshat} $  the subcategory of
$ \mathcal{T}_{\otimeshat} $  whose objects are all the Hopf
algebras in  $ \mathcal{T}_{\otimeshat} $  and whose morphisms
are all the Hopf algebra morphisms in $ \mathcal{T}_{\otimeshat} $.
                                            \par \noindent
 Second, let  $ \mathcal{P}_{\otimestilde} $  be the category whose objects are all topological  $ \kh $--modules  isomorphic to modules of the type  $ {\kh}^E $  (with the product topology) for some set 
$ E \, $,  and whose morphisms are the  $ \kh $--linear  continuous maps.  Again, this is a tensor category w.r.t.~the tensor product  $ \, P_1
\otimestilde P_2 \, $  defined as the comple\-tion of the algebraic
tensor product $ \, P_1 \otimes_{\kh} P_2 \, $  w.r.t.~the weak topology:
thus  $ \, P_i \cong {\kh}^{E_i} $  ($ i = 1 $,  $ 2 $)  yields  $ \, P_1
\otimestilde P_2 \cong {\kh}^{E_1 \times E_2} \, $  (for  $ P_1 $,  $ P_2
\in \mathcal{P}_{\otimestilde} $).  We call  $ \HA_{\otimestilde} $  the
subcategory of  $ \mathcal{P}_{\otimestilde} $  with objects the Hopf
algebras in  $ \mathcal{P}_{\otimestilde} \, $,  and with morphisms
the Hopf algebra morphisms in  $ \mathcal{P}_{\otimestilde} $.
                                         \par \noindent
A  {\sl quantum group\/}  for us will be a Hopf algebra, in either of
$ \HA_{\otimeshat} $  or  $ \HA_{\otimestilde} \, $,  having a special semiclassical limit.  The exact definition is the following:

\begin{definition}  (cf.~[7], \S~7)
                                         \par
   {\it (a)}  We call QUEA any  $ \, H \in \HA_{\otimeshat} $
such that  $ \, H_0 := H \big/ \h H \, $  is a co-Poisson Hopf algebra
isomorphic to  $ U(\gerg) $  for some finite dimensional Lie bialgebra
$ \gerg $  (over  $ \Bbbk $);  then we write  $ \, H = U_\h(\gerg) \, $,
\, and say  $ H $  is a  {\sl quantization\/}  of  $ U(\gerg) $.  We call
$ \QUEA $  the full tensor subcategory of  $ \, \HA_{\otimeshat} $  whose
objects are QUEA, relative to all possible  $ \gerg \, $.
                                         \par
   {\it (b)}  We call QFSHA any  $ \, K \in \HA_{\otimestilde} $
such that  $ \, K_0 := K \big/ \h K \, $  is a topological Poisson
Hopf algebra isomorphic to  $ F[[G]] $  for some finite dimensional
formal Poisson group  $ G $  (over  $ \Bbbk $);  then we write  $ \,
H = F_\h[[G]] \, $,  \, and say  $ K $  is a  {\sl quantization\/}
of  $ F[[G]] $.  We call  $ \QFSHA $  the full tensor subcategory
of  $ \, \HA_{\otimestilde} $  whose objects are QFSHA, relative
to all possible  $ G \, $.
\end{definition}    

\subsection{Remarks}    
 If  $ \, H \in \HA_{\otimeshat} $  is such that its semiclassical
limit  $ \, H_0 := H \big/ \h H \, $  as a Hopf algebra is isomorphic
to  $ U(\gerg) $  for some Lie algebra  $ \gerg \, $,  then  $ \, H_0
= U(\gerg) \, $  is also a  {\sl co-Poisson\/}  Hopf algebra w.r.t.~the
Poisson cobracket  $ \delta $  defined as follows: if  $ \, x \in H_0
\, $  and  $ \, x' \in H \, $  gives  $ \, x = x' + \h \, H \, $, 
\, then  $ \, \delta(x) := \big( \h^{-1} \, \big( \Delta(x') -
\Delta^{\text{op}}(x') \big) \big) + \h \, H \otimeshat H \, $. 
Thus, in particular   --- by Theorem 2 in \S 3 of Drinfeld's paper
[6] ---   the restriction of  $ \delta $  makes  $ \gerg $  into a
Lie bialgebra.  In  Definition 3.1{\it (a)},  the co-Poisson structure
considered on  $ H_0 $  is nothing but the one arising in this way. 
Similarly, if  $ \, K \in \HA_{\otimestilde} $  is such that its
semiclassical limit  $ \, K_0 := K \big/ \h K \, $  is a topological
Poisson Hopf algebra isomorphic to  $ F[[G]] $  for some formal group 
$ G $  then  $ \, K_0 = F[[G]] \, $  is also a topological  {\sl
Poisson\/}  Hopf algebra w.r.t.~the Poisson bracket  $ \{\,\ ,\ \} $ 
defined as follows: if  $ \, x $,  $ y \in K_0 \, $  and  $ \, x' $, 
$ y' \in K \, $  give  $ \, x = x' + \h \, K \, $  and  $ \, y = y' +
\h \, K $, then  $ \, \{x,y\} := \big( \h^{-1} (x' \, y' - y' \, x')
\big) + \h \, K \, $.  Then, in parti\-cular,  $ G $  is a  {\sl
Poisson\/}  formal group.  And again, in  Definition 3.1{\it (b)}, 
the Pois\-son structure considered on  $ K_0 $  is exactly   
     \hbox{the one that arises from this construction.}

\subsection{Drinfeld's functors}   
 Let  $ H $  be a (topological)
Hopf algebra over  $ \kh $.  Letting  $ \, J_{\scriptscriptstyle H}
:= \text{\it Ker}\,(\epsilon_{\scriptscriptstyle H}) \, $  and  $ \, I_{\scriptscriptstyle H} := \epsilon_{\scriptscriptstyle H}^{-1} \big(
\h \, \kh \big) = J_{\scriptscriptstyle H} + \h \, H \, $,  we set
  $$  H^\times \, := \, {\textstyle \sum\nolimits_{n \geq 0}}
\h^{-n} {I_{\scriptscriptstyle H}}^{\!n} = {\textstyle
\sum\nolimits_{n \geq 0}} {\big( \h^{-1} I_{\scriptscriptstyle H}
\big)}^n = {\textstyle \bigcup\limits_{n \geq 0}}
{\big( \h^{-1} I_{\scriptscriptstyle H} \big)}^n =
{\textstyle \sum\nolimits_{n \geq 0}} \h^{-n}
{J_{\scriptscriptstyle H}}^{\!n}  $$
 \vskip-5pt
\noindent
 which is a subspace of  $ \, \Bbbk((\h)) \otimes_{\kh} H \, $. 
Then we define
  $$  H^\vee :=  \h\text{--adic completion of the
$ \kh $--module }  H^\times  \; .  $$
 On the other hand, for each  $ \, n \in \N_+ \, $,  define  $ \;
\Delta^n \colon H \longrightarrow H^{\otimes n} \; $  by  $ \,
\Delta^1 := \id_{\scriptscriptstyle H} $  and  $ \, \Delta^n :=
\big( \Delta \otimes \id_{\scriptscriptstyle H}^{\,\otimes (n-2)}
\big) \circ \Delta^{n-1} \, $  if  $ \, n \geq 2 \, $,  and set
$ \, \delta_0 := \delta_\emptyset \, $,  and  $ \, \delta_n :=
{(\id_{\scriptscriptstyle H} - \epsilon)}^{\otimes n} \circ
\Delta^n \, $,  for all  $ \, n \in \N_+ \, $.  Then we define
  $$  H' := \big\{\, a \in H \,\big\vert\; \delta_n(a) \in
h^n H^{\otimes n} \;\; \forall\, n \in \N \,\big\}  \qquad
\big( \subseteq H \, \big ) \; .   $$
 \vskip3pt
\noindent
 Note that the definition of  $ H^\vee $  is pretty direct.  In particular,
we specify how it can be generated (topologically), namely it is the
(complete topological) unital  $ \kh $--subalgebra  of  $ \, \Bbbk((\h))
\otimes_{\kh} H \, $  generated by  $ \h^{-1} J_{\scriptscriptstyle H} $
or  $ \h^{-1} I_{\scriptscriptstyle H} \, $.  In contrast, the definition
of  $ H' $  is quite implicit: roughly speaking, it is the set of solution of a system with countably many equations
(specified in terms of  $ \h $--adic  valuation).  Nevertheless, the two definitions are strictly related, in a sense made explicit by Proposition 2.6 below.
                                            \par \noindent
 Now we state the Quantum Duality Principle (=QDP) for quantum groups:

\begin{theorem} 
 ({\sl see [6], and [11] for a proof)}  {\it The
assignments  $ \, H \mapsto H^\vee \, $  and  $ \, H \mapsto H' \, $,
respectively, define tensor functors  $ \, \QFSHA \longrightarrow \QUEA
\, $  and  $ \, \QUEA \longrightarrow \QFSHA \, $,  which are inverse
to each other.  Indeed,  for all  $ \, U_\h(\gerg) \in \QUEA \, $  and
all  $ \, F_\h[[G]] \in \QFSHA \, $  one has
  $$  {U_\h(\gerg)}' \Big/ \h \, {U_\h(\gerg)}' = F[[G^*]]
\hskip37pt  {F_\h[[G]]}^\vee \Big/ \h \, {F_\h[[G]]}^\vee
= U(\gerg^*)  $$
that is, if  $ \, U_\h(\gerg) $  is a quantization of  $ \, U(\gerg) $
then  $ \, {U_\h(\gerg)}' $  is one of  $ \, F[[G^*]] $,  and if  $ \,
F_\h[[G]] $  is a quantization of  $ \, F[[G]] $  then  $ \,
{F_\h[[G^*]]}^\vee $  is one of  $ \, U(\gerg^*) \, $.}
\end{theorem}    

\noindent
 In addition, Drinfeld's functors respect Hopf duality, in the following
sense:

\begin{proposition}    
 {\sl (see [11], Proposition 2.2)}  {\it Let
$ \, U_\h \in \QUEA \, $,  $ \, F_\h \in \QFSHA $  and let  $ \, \pi
\, \colon \, U_\h \times F_\h \!\longrightarrow \kh \, $  be a perfect
Hopf pairing whose specialization at  $ \, \h = 0 \, $  is perfect as well.
Then  $ \pi $  induces   --- by restriction on l.h.s.~and scalar extension
on r.h.s.~---   a perfect Hopf pairing  $ \, {U_\h}' \times {F_\h}^{\!\vee}
\!\longrightarrow \kh \, $  whose specialization at  $ \h = 0 $  is again
perfect.}
\end{proposition}    

\noindent
 In other words, the above result ensures that, if one starts from a pair
made by a QUEA and a QFSHA which are dual to each other, and then applies
Drinfeld's functors to both terms of the pair, then one obtains another
pair   --- now with QFSHA first and QUEA second ---   with the same
property.  In this sense, the two Drinfeld's functors are ``dual to
each other''.

\subsection{Quantum subgroups and quantum homogeneous spaces}   
 From now on, let  $ G $  be a formal Poisson group,  $ \, \gerg :=
\text{\it Lie}(G) $  its tangent Lie bialgebra.  We assume a quantization
of  $ G $  is given, in the sense that a QFSHA  $ \fhg $  quantizing 
$ F[[G]] $  and a QUEA  $ \uhg $  quantizing  $ U(\gerg) $  are given
such that, in addition,  $ \, \fhg \cong {\uhg}^* := \text{\sl Hom}_{\,
\kh}\big(\uhg,\kh\big) \, $  as topological Hopf algebras; the latter
requirement is equivalent to fixing a perfect Hopf algebra pairing between 
$ \fhg $  and  $ \uhg \, $  whose specialization at  $ \, \h = 0 \, $ 
be perfect too.  This assumption is not restrictive, because [8]
such a  $ \uhg $  always exists, and then one can take  $ \, \fhg
:= {\uhg}^* \, $.  We denote by  $ \, \pi_{F_\h} \colon \fhg
\loongtwoheadrightarrow F[[G]] \, $   and  $ \, \pi_{U_\h} \colon
\uhg \loongtwoheadrightarrow U(\gerg) \, $  the specialization maps,
and we set  $ \, F_\h := \fhg \, $,  $ \, U_\h := \uhg \, $.
                                               \par \noindent
   Let  $ K $  be a formal subgroup of  $ G \, $,  and  $ \, \gerk :=
\text{\it Lie}(K) \, $.  As quantization of  $ K $  and/or of  $ \, G
\big/ K \, $,  we mean a quantization of any one of the four algebraic
objects  $ \calI $,  $ \calC $,  $ \gerI $  and  $ \gerC $  associated
to them in Subsection 2.2, that is either of
 \vskip-13pt
  $$  \hbox{ \hskip-9pt   $ \begin{matrix}
   (a)  &  \, \text{\ a left ideal, coideal \ }  \calI_\h \ideal_\ell
\coideal \; F_\h[[G]]  \, \text{\ such that \ }  \hfill  \\
        &  \qquad \phantom{\Big|_|} \qquad  \calI_\h \big/ \h \,
\calI_\h \, \cong \, \pi_{F_\h}(\calI_\h) \, = \, \calI   \hfill  \\
   (b)  &  \, \text{\ a subalgebra, left coideal \ }  \calC_\h \leq^1
\coideal_\ell \, F_\h[[G]]  \, \text{\ such that \ }  \hfill  \\
        &  \qquad \phantom{\Big|_|} \qquad \calC_\h \big/ \h \,
\calC_\h \, \cong \, \pi_{F_\h}(\calC_\h) \, = \, \calC   \hfill  \\
   (c)  &  \, \text{\ a left ideal, coideal \ } \gerI_\h \ideal_\ell
\coideal \;\; U_\h(\gerg)  \, \text{\ such that\ }   \hfill  \\
        &  \qquad \phantom{\Big|_|} \qquad  \gerI_\h \big/ \h \,
\gerI_\h \, \cong \, \pi_{U_\h}(\gerI_\h) \, = \, \gerI   \hfill  \\
   (d)  &  \, \text{\ a subalgebra, left coideal \ } \gerC_\h \leq^1
\coideal_\ell \, U_\h(\gerg)  \, \text{\ such that\ }   \hfill  \\
        &  \qquad \phantom{\Big|_|} \qquad  \gerC_\h \big/ \h \,
\gerC_\h \, \cong \, \pi_{U_\h}(\gerC_\h) \, = \, \gerC   \hfill
    \end{matrix} $ }   \eqno (3.1)  $$
 \vskip-9pt
\noindent
 In (3.1) above, the constraint  $ \; \calI_\h \big/ \h \, \calI_\h \, \cong \, \pi_{F_\h}(\calI_\h) \, = \, \calI \; $  means the following.
By construction  $ \; \calI_\h \! \longrightarrow \! \fhg
\,{\buildrel{\pi_{F_\h}}\over{\loongtwoheadrightarrow}}\, \fhg
\Big/ \h \, \fhg \, \cong \, F[[G]] \, $,  and the composed map 
$ \, \calI_\h \! \longrightarrow \! F[[G]] $  factors through  $ \, \calI_\h \big/ \h \, \calI_\h \, $;  \, then we ask that the induced
map  $ \, \calI_\h \big/ \h \, \calI_\h \longrightarrow F[[G]] \, $  be
a bijection onto  $ \pi_{F_\h}(\calI_\h) \, $,  and that the latter
do coincide with  $ \calI \, $;  of course this bijection will also
respects all Hopf operations, because  $ \pi_{F_\h} $  does.
Similarly for the other conditions.
                                            \par \noindent
 Moreover, let
   $ X \! \in \! \{ \calI, \calC,
\gerI, \gerC \} $,  $ S_\h \! \in \! \big\{ \fhg, \uhg \!\big\} $.
Since $ \, \pi_{S_\h}(X_\h) = X_\h \Big/ \! \big( X_\h \cap \h \,
S_\h \big) \, $,  the property  $ \, X_\h \Big/ \h \, X_\h \cong
\pi_{S_\h}(X_\h) = X \, $  is equivalent to  $ \, X_\h \cap \h \, S_\h 
= \h \, X_\h \, $.
So our quantum objects can also be characterized by
 \vskip-13pt
  $$  \hskip-1pt   \hbox{ $ \begin{matrix}
   &  (a) \quad \hskip2pt  \calI_\h \ideal_\ell \coideal \; \fhg
\hfill  &  \,\quad\;  \calI_\h \cap \h \, \fhg \, = \, \h \, \calI_\h
&  \,\quad\;  \calI_\h \big/ \h \, \calI_\h \, = \, \calI  \\
   &  (b) \quad \hskip2pt  \calC_\h \leq^1 \coideal_\ell \, \fhg
\hfill  &  \,\quad\;  \calC_\h \cap \h \, \fhg \, = \, \h \, \calC_\h
&  \,\quad\;  \calC_\h \big/ \h \, \calC_\h \, = \, \calC  \\
   &  (c) \quad \hskip2pt  \gerI_\h \ideal_\ell \coideal \;\; U_\h(\gerg)  \hfill  &  \,\quad\;  \gerI_\h \cap \h \, U_\h(\gerg) \, = \, \h \, \gerI_\h  &  \,\quad\;  \gerI_\h \big/ \h \, \gerI_\h \, = \, \gerI  \\
   &  (d) \quad \hskip2pt  \gerC_\h \leq^1 \coideal_\ell \, U_\h(\gerg)   \hfill  &  \,\quad\;  \gerC_\h \cap \h \, U_\h(\gerg) \, = \, \h \, \gerC_\h  &  \,\quad\;  \gerC_\h \big/ \h \, \gerC_\h \, = \, \gerC
               \end{matrix} $ }   \eqno (3.1)'  $$
 \vskip1pt
\noindent
instead of (3.1).  Note that  $ \, \calI = \calI(K) \, $  and  $ \, \gerC = \gerC(K) \, $  provide an ``alge\-braization'' (in global and local terms respectively) of the subgroup  $ K $,  more than of the homogeneous space  $ G \big/ K \, $; conversely,  $ \, \calC = \calC(K) \, \big( \cong F \big[ \big[ G/K \big] \big] \, \big) \, $  and  $ \, \gerI = \gerI(K) \, $  provide an ``algebraization'' (of global and local type
respectively) of  $ G \big/ K  $,  more than of  $ K \, $.  For this reason, in the sequel we shall loosely refer to  $ \calI_\h $  and 
$ \gerC_\h $ as to  {\sl ``quantum (formal) subgroups''},  and to 
$ \calC_\h $  and  $ \gerI_\h $  instead as to  {\sl ``quantum
(formal) homogeneous spaces''}.

\vskip5pt

\noindent
 One could ask whether quantum subgroups and quantum homogeneous spaces do exist.  Actually, from the
very definitions one immediately finds a square  {\sl necessary condition}.  In fact, next Lemma proves that the (formal) subgroup
of  $ G $  obtained as specialization of a quantum (formal) subgroup is  {\sl coisotropic},  and a (formal) homogeneous  $ G $--space  obtained as specialization of a quantum (formal) homogeneous space is a  {\sl Poisson quotient}.  This is quite a direct generalization of the situation for quantum groups, where one has that specializing a quantum group always gives a  {\sl Poisson\/}  group.

\begin{lemma}   
 {\it Let  $ K $  be a formal subgroup of
$ G $,  and assume a quantization  $ \, \calI_\h $,  $ \calC_\h $,
$ \gerI_\h $  or  $ \, \gerC_\h $  of  $ \, \calI $,  $ \calC $,
$ \gerI $  or  $ \, \gerC $  respectively be given as above.  Then
the subgroup  $ K $  is coisotropic, and the  $ G $--space
$ \, G \big/ K $  is a Poisson quotient.}
\end{lemma}   

\noindent
 Note that, at the quantum level, one looses either commutativity
or cocommutativity; then, one-sided ideals (or coideals) are not
automatically two-sided!  This enters in the definitions above,
in that we require some objects to be  {\sl one\/}-sided ideals
(coideals)   --- taking  {\it left\/}  rather than  {\it right\/} 
ones is just a matter of choice.  If one takes  {\sl two\/}-sided 
ones instead, the like of Lemma 3.1 is that  $ K $  be a  {\sl
Poisson subgroup\/}  (for  $ \calI(K) $  is a Poisson ideal).   

\subsection{The existence problem}      
 The  {\sl existence\/}  of any of the four possible objects providing
a quantization of a coisotropic subgroup (or of the associated Poisson
quotient) is an open problem.  Etingof and Kahzdan [9] gave a positive
answer for the subclass of those  $ K $  which are also  {\sl Poisson
subgroups\/}  (which amounts to  $ \, \gerk := \text{\it Lie}(K) \, $ 
being a Lie subbialgebra).  Many other examples of quantizations exist
too.
                                                   \par \noindent
 Yet, the four existence problems are equivalent: i.e., as one solves
any one of them, a solution follows for the remaining ones.  Indeed,
one has:

 \vskip4pt

   {\it ---  $ \underline{\hbox{{\it (a)}  $ \Longleftrightarrow $
{\it (d)}  \, and \,  {\it (b)}  $ \Longleftrightarrow $  {\it (c)}}} $:}
\, if  $ \calI_\h $  exists as in  {\it (a)},  then  $ \, \gerC_\h :=
{\calI_\h}^{\!\perp} \, $  (hereafter orthogonality is meant w.r.t.~the fixed Hopf pairing between  $ \fhg $  and  $ U_\h(\gerg) \, $)  enjoys
the properties in  {\it (d)\/};  \, conversely, if  $ \gerC_\h $  exists
as in  {\it (d)},  then  $ \, \calI_\h := {\gerC_\h}^{\!\perp} \, $  enjoys the properties in {\it (a)}.  Similarly, the equivalence \,
{\it (b)}  $ \Longleftrightarrow $  {\it (c)} \; follows from a like orthogonality argument.

 \vskip2pt

   {\it ---  $ \underline{\hbox{{\it (a)}  $ \Longleftrightarrow $
{\it (b)}  \, and \,  {\it (c)}  $ \Longleftrightarrow $  {\it (d)}}} $:}
\, if  $ \calI_\h $  exists as in  {\it (a)},  then  $ \, \calC_\h :=
\calI_\h^{\,\text{\it co}\calI_\h} \, $  is an object like in  {\it
(b)\/};  \, on the other hand, if  $ \calC_\h $  as in  {\it (b)\/}
is given, then  $ \, \calI_\h := \fhg \cdot \calC_\h^{\,+} \, $
enjoys all properties in  {\it (a)\/}  (notation of Subsection 2.2). 
The equivalence \,  {\it (c)}  $ \! \Longleftrightarrow \! $  {\it (d)}
\; stems from a like argument.

\subsection{Basic assumptions}     
 Hereafter  {\sl we assume that quantizations  $ \calI_\h \, $,  $ \calC_\h \, $,  $ \gerI_\h \, $,  $ \gerI_\h $  as in (3.1) are given and that they be linked by relations like  {\it (1)--(2)\/}  in Subsection 2.2, i.e.}
  $$  \hskip-8pt   \hbox{ $ \begin{matrix}
   &  \hbox{\it (i)}  \;\quad  \calI_\h = {\gerC_\h
\phantom{)}}^{\!\!\!\!\perp} \, ,  \;\quad  \gerC_\h
= {\calI_\h \phantom{)}}^{\!\!\!\perp}  \quad
  &  \hbox{\it (ii)}  \;\quad  \gerI_\h = {\calC_\h
\phantom{)}}^{\!\!\!\perp} \, ,  \;\quad  \calC_\h =
{\gerI_\h \phantom{)}}^{\!\!\!\!\perp}  \\
   &  \hbox{\it (iii)}  \;\;  \calI_\h = F_\h \cdot \calC_\h^{\,+} \, ,
\;\;  \calC_\h = {F_\h}^{\!\text{\it co}\calI_\h}  \quad
  &  \hbox{\it (iv)}  \;\;  \gerI_\h = U_\h \cdot \gerC_\h^{\,+} \, ,
\;\;  \gerC_\h = {U_\h}^{\!\text{\it co}\gerI_\h}  \\
   \end{matrix} $ }   \hskip6pt (3.2)  $$
 \vskip2pt
\noindent
 In fact, one of the objects is enough to have all others, in
such a way that the previous assumption holds.  Indeed, if  $ \,
\hbox{\it co}\hskip1pt\calS \, $  is the set of coisotropic subgroup of  $ G $,  let  $ \, Y_\h \big( \hbox{\it co}\hskip1pt\calS \big)
:= \big\{ Y_\h(\gerk) \big\}_{\gerk \in \, co\calS} \, $  for all 
$ \, Y \in \big\{ \calI, \calC, \gerI, \gerC \big\} \, $.  The
four equivalences  {\it (a)}  $ \Longleftrightarrow $  {\it (d)},
{\it (b)}  $ \Longleftrightarrow $  {\it (c)},  {\it (a)}
$ \Longleftrightarrow $  {\it (b)}  \, and  {\it (c)}
$ \Longleftrightarrow $  {\it (d)}  \, above are given
by bijections  $ \, \calI_\h \big( \hbox{\it
co}\hskip1pt\calS \big) \longleftrightarrow \gerC_\h
\big( \hbox{\it co} \hskip1pt\calS \big) \, $,  $ \,
\calC_\h \big( \hbox{\it co}\hskip1pt\calS \big)
\longleftrightarrow \gerI_\h \big( \hbox{\it co}
\hskip1pt\calS \big) \, $,  $ \, \calI_\h \big(
\hbox{\it co}\hskip1pt\calS \big) \longleftrightarrow
\calC_\h \big( \hbox{\it co}\hskip1pt\calS \big) \, $
and  $ \, \gerI_\h \big( \hbox{\it co}\hskip1pt\calS
\big) \longleftrightarrow \gerC_\h \big( \hbox{\it co}
\hskip1pt\calS \big) \, $  which form a  {\sl commutative\/}  square.   
In fact, each of these maps, or their inverse, is of type  $ \, X_\h \mapsto X_\h^\perp \, $,  $ \, A_\h \mapsto H_\h A_\h^+ \, $  or  $ \, K_\h \mapsto H_\h^{\text{\it co}K_\h} \, $  (see Subsection 2.2): since  $ \, X_\h \subseteq \big(X_\h^\perp\big)^\perp \, $  and  $ \, A_\h \subseteq H_\h^{\text{\it co}(H_\h A_\h^+)} \, $  in general, and these inclusions are identities at  $ \, \h = 0 \, $,  one gets  $ \, X_\h = \big( X_\h^\perp \big)^\perp \, $  and  $ \, A_\h = H_\h^{\text{\it co}(H_\h A_\h^+)} \, $.   
                                          \par \noindent
   Note that  $ \calI_\h \big( \hbox{\it co}\hskip1pt
\calS \big) $,  $ \gerC_\h \big( \hbox{\it co}\hskip1pt\calS
\big) $,  $ \calC_\h \big( \hbox{\it co}\hskip1pt\calS \big) $
and  $ \gerI_\h \big( \hbox{\it co}\hskip1pt\calS \big) $  are
again lattices with respect to inclusion, hence they will 
be thought of as categories too.

\subsection{Remark}     
 If a quadruple  $ \, \big( \calI_\h \, , \,
 \calC_\h \, , \, \gerI_\h \, , \, \gerC_\h \big) \, $  is given which enjoys all properties in the first and the second column of  $ (3.1)' $,
along with relations (3.2), then one easily checks that the four specialized objects  $ \, \calI := \calI_\h\big|_{\h=0} \, $, $ \,
\calC := \calC_\h\big|_{\h=0} \, $,  $ \, \gerI := \gerI_\h\big|_{\h=0}
\, $  and  $ \, \gerC := \gerC_\h\big|_{\h=0} \, $  verify relations
{\it (1)\/}  and  {\it (2)\/}  in Subsection 2.2.  Therefore, these four objects define just  {\sl one\/}  single pair  {\sl (coisotropic subgroup, Poisson quotient)},  and the quadruple  $ \, \big( \calI_\h \, , \, \calC_\h \, , \, \gerI_\h \, , \, \gerC_\h \big) \, $  yields a quantization of the latter in the sense of Subsection 3.4.  
$ \quad \diamondsuit $

\subsection{General program}     
 From the setup of Subsection 2.2, we follow this scheme: 
 \vskip-17pt
  $$  \begin{matrix}
   \hbox{\it (a)}  \phantom{\bigg|} \quad  &  \big( F[[G]] \supseteq \big)
\  &  \calI  &  \;{\buildrel (1) \over \llongrightarrow}\;  &  \calI_\h  &  \;{\buildrel (2) \over \llongrightarrow}\;  &  {\calI_\h}^{\!\curlyvee}  &  \,{\buildrel (3) \over \llongrightarrow}\;  &  {\calI_0}^{\!\curlyvee}  &
\!\ \big( \subseteq U\big(\gerg^*\big) \,\big)  \\
   \hbox{\it (b)}  \phantom{\bigg|} \quad  &  \big( F[[G]] \supseteq \big)
\  &  \calC  &  \;{\buildrel (1) \over \llongrightarrow}\;  &  \calC_\h
&  \;{\buildrel (2) \over \llongrightarrow}\;  &  {\calC_\h}^{\!\!\triangledown}  &  \,{\buildrel (3) \over \llongrightarrow}\;  &  {\calC_0}^{\!\!\triangledown}  &  \!\ \big( \subseteq U\big(\gerg^*\big) \,\big)  \\
   \hbox{\it (c)}  \phantom{\bigg|} \quad  &  \big( U(\gerg) \supseteq \big) \  &  \gerI  &  \;{\buildrel (1) \over \llongrightarrow}\;  &  \gerI_\h  &  \;{\buildrel (2) \over \llongrightarrow}\;  &  {\gerI_\h}^{\! !}  &  \,{\buildrel (3) \over \llongrightarrow}\;  &  {\gerI_0}^{\! !}  &
\!\ \big( \subseteq F[[G^*]] \,\big)  \\
   \hbox{\it (d)}  \phantom{\bigg|} \quad  &  \big( U(\gerg) \supseteq \big) \  &  \gerC  &  \;{\buildrel (1) \over \llongrightarrow}\;  &  \gerC_\h  &  \;{\buildrel (2) \over \llongrightarrow}\;  &  {\gerC_\h}^{\!\!\Lsh}  &  \,{\buildrel (3) \over \llongrightarrow}\;  &  {\gerC_0}^{\!\!\Lsh}  &
\!\ \big( \subseteq F[[G^*]] \,\big)
   \end{matrix}  $$

\vskip-6pt

\noindent
Its meaning is the following.  Starting from the first column (in left-hand side) we move a first step   --- arrows (1) ---   which is some quantization process.  Instead, the last step   --- arrows (3) ---   is a specialization (at  $ \, \h = 0 \, $)  process. 
In between, the middle step   --- arrows (2) ---   demands some new idea: here we shall apply some suitable ``adaptations'' of Drinfeld's functors to the quantizations (of a coisotropic subgroup or a Poisson quotient) obtained from step (1).  Roughly, the idea is to take the suitable Drinfeld's functor on the quantum group   ---  $ \fhg $  or 
$ \uhg $  ---   and to ``restrict it'' to the quantum sub-object
given by step (1).  The points to show then are
 \vskip2pt
   {\sl  $ \underline{\text{First}} $:} \, each one of the four objects in the third column above   --- that is, provided by arrows (2) ---   is one of the four algebraic objects which describe a quantum (closed formal) coisotropic subgroup or Poisson quotient of  $ G^* \, $.  Namely, the correspondence of ``types'' is
 \vskip1pt
   \centerline{ $
\text{\it (a)} \, =\joinrel\Longrightarrow \, \text{\it (c)}
\;\; ,  \hskip17pt
\text{\it (b)} \, =\joinrel\Longrightarrow \, \text{\it (d)}
\;\; ,  \hskip17pt
\text{\it (c)} \, =\joinrel\Longrightarrow \, \text{\it (a)}
\;\; ,  \hskip17pt
\text{\it (d)} \, =\joinrel\Longrightarrow \, \text{\it (b)} $ }
 \vskip2pt
\noindent   
 where notation  $ \; \text{\it (x)} =\joinrel\Longrightarrow
\! \text{\it (y)} \; $   means that an object of ``type  {\it (x)\/}''   --- referring to the classification of (3.1) or of  $ (3.1)'$  ---   yields an object of ``type  {\it (y)\/}''.
 \vskip2pt
   {\sl  $ \underline{\text{Second}} $:} \, the four formal subgroups or
homogeneous spaces of  $ G^* $  obtained above provide only  {\sl one\/}  single pair  {\sl (subgroup, homogeneous space)}.
 \vskip2pt
   {\sl  $ \underline{\text{Third}} $:} \, if we start from a coisotropic subgroup  $ K \, $,  and/or a Poisson quotient  $ \, G \big/ K \, $,  of  $ \, G \, $,  then the Poisson quotient and/or the subgroup of  $ \, G^* $  obtained above are  $ \, G^* \big/ K^\perp
\, $  and/or  $ K^\perp $  (cf.~Definition 2.1) respectively.

 \vskip1,0truecm

\section{Drinfeld-like functors on quantum subgroups
and quantum Poisson quotients}

\vskip5pt

\subsection{Restricting Drinfeld's functors}     
 The main idea in the program sketched in Subsection 3.8 is the intermediate step   --- provided by arrows (2)   --- namely that of ``restricting'' Drinfeld's functors, originally defined for quantum groups, to the quantum sub-objects we are interested in.  To this end, the ``right'' definition is the following:

\begin{definition}     
Let  $ \, J := \text{\it Ker}\,(\epsilon_{\scriptscriptstyle F_\h[[G]]}) \, $, 
$ \, I := J + \h \, F_\h[[G]] \; $.
 \vskip2pt
\noindent
 {\it (a)} \hskip7pt  $ \displaystyle{ {\calI_\h}^{\!\curlyvee}
\; := \;  {\textstyle \sum\limits_{n=1}^\infty} \, \h^{-n} \cdot I^{n-1}
\cdot \calI_\h  \; = \;  {\textstyle \sum\limits_{n=1}^\infty} \, \h^{-n}
\cdot J^{n-1} \cdot \calI_\h } $
 \vskip1pt
\noindent
 {\it (b)} \hskip7pt  $ \displaystyle{ {\calC_\h}^{\!\! \triangledown}
\; := \;  \calC_\h \, + \, {\textstyle \sum\limits_{n=1}^\infty} \,
\h^{-n} \cdot {\big( \calC_\h \cap I \,\big)}^n \, =  \;  \kh \cdot 1
\, + \, {\textstyle \sum\limits_{n=1}^\infty} \, \h^{-n} \cdot {\big(
\calC_\h \cap J \,\big)}^n } $
 \vskip1pt
\noindent
 {\it (c)} \hskip5pt  $ \displaystyle{ {\gerI_\h}^{\! !}  \, :=
\,  \bigg\{\, x \in \gerI_\h \,\bigg\vert\, \delta_n(x) \in \h^n
{\textstyle \sum\limits_{s=1}^n} \, {U_\h}^{\!\otimeshat (s-1)}
\otimeshat \gerI_\h \otimeshat {U_\h}^{\!\otimeshat (n-s)} \! , \,
\forall\; n \in \N_+ \bigg\} } $
 \vskip1pt
\noindent
 {\it (d)} \hskip7pt  $ \displaystyle{ {\gerC_\h}^{\!\!\Lsh} \;
:= \;  \Big\{\, x \in \gerC_\h \;\Big\vert\; \delta_n(x) \in \h^n
\, {U_\h}^{\!\otimeshat (n-1)} \otimeshat \gerC_\h \, , \; \forall\;
n \in \N_+ \,\Big\} } $
\end{definition}   

\noindent
Indeed, directly by definitions one has that
 \vskip-16pt
  $$  {\calI_\h}^{\!\curlyvee} \supseteq \calI_\h \;\;\; ,
\hskip25pt  {\calC_\h}^{\!\!\triangledown} \supseteq \calC_\h
\;\;\; ,  \hskip25pt  {\gerI_\h}^{\! !} \subseteq \gerI_\h \;\;\; ,
\hskip25pt  {\gerC_\h}^{\!\!\Lsh} \subseteq \gerC_\h \;\;\; .  $$
 \vskip-4pt
\noindent   
 But even more, a careful (yet easy) analysis of definitions and of the
relationship between each quantum groups and its relevant sub-objects
shows that, in force of  $ (3.1)' $   --- in particular, the mid column
there ---   one has
 \vskip-16pt
  $$  \calI_\h \, = \, {\calI_\h}^{\! \curlyvee} \cap F_\h \; ,  \hskip13pt
\calC_\h \, = \, {\calC_\h}^{\!\! \triangledown} \cap F_\h \; ,  \hskip13pt
{\gerI_\h}^{\! !} = \, \gerI_\h \cap {U_\h}' \; ,  \hskip13pt
{\gerC_\h}^{\!\!\Lsh} = \, \gerC_\h \cap {U_\h}'  $$
 \vskip-4pt
\noindent
 This proves that, in very precise sense, Definition 4.1 really provides a
``restriction'' of Drinfeld's functors from quantum groups to our quantum
subgroups   --- in cases  {\it (a)\/}  and  {\it (d)} ---   or quantum
Poisson quotients   --- in cases  {\it (b)\/}  and  {\it (c)}.  This
also motivates the notation: indeed, the symbols  $ {}^\curlyvee $  and
$ {}^\triangledown $  are (or should be) remindful of  $ {}^\vee $,
while  $ {}^! $  and  $ {}^\Lsh $  are remindful of  $ {}' $.

\vskip3pt

\noindent
 We can now state the QDP for coisotropic subgroups and
Poisson quotients:

\vskip13pt

\begin{theorem}     
 ({\sl ``QDP for Coisotropic Subgroups and
Poisson Quotients''})
 \vskip1pt
   {\it (a) \, Definition 4.1 provides category equivalences
 \vskip-14pt
  $$  \displaylines{
   {(\ )}^\curlyvee \colon \, \calI_\h \big( \hbox{\it co}
\hskip1pt\calS(G) \!\big) \, {\buildrel \cong \over \longrightarrow}
\, \gerI_\h \big( \hbox{\it co}\hskip1pt\calS(G^*) \!\big) \; ,
\quad {(\ )}^\triangledown \colon \, \calC_\h \big( \hbox{\it
co}\hskip1pt\calS(G) \!\big) \, {\buildrel \cong \over \longrightarrow}
\, \gerC_\h \big( \hbox{\it co}\hskip1pt\calS(G^*) \!\big)  \cr
   {(\ )}^! \, \colon \, \gerI_\h \big( \hbox{\it co}\hskip1pt
\calS(G) \!\big) \, {\buildrel \cong \over \longrightarrow} \,
\calI_\h \big( \hbox{\it co}\hskip1pt\calS(G^*) \!\big) \; ,
\quad  {(\ )}^\Lsh \colon \, \gerC_\h \big( \hbox{\it co}
\hskip1pt\calS(G) \!\big) \, {\buildrel \cong \over \longrightarrow}
\, \calC_\h \big( \hbox{\it co}\hskip1pt \calS(G^*) \!\big)  \cr }  $$
 \vskip-3pt
\noindent
 along with the similar ones with  $ G $  and  $ G^* $  interchanged,
such that  $ \, {(\ )}^! $  and  $ \, {(\ )}^\curlyvee $  are
inverse to each other, and  $ \, {(\ )}^\Lsh $  and  $ \,
{(\ )}^\triangledown $  are inverse to each other.
 \vskip4pt
   (b) \,  {\sl (the QDP)}  For any  $ \, K \in
\hbox{\it co}\hskip1pt\calS(G) \, $  and  $ \,
\gerk := \text{\it Lie}(K) \, $,  we have
 \vskip-13pt
  $$  \begin{matrix}
   \calI(K)_\h^{\,\curlyvee} \mod \h \, {F_\h[[G]]}^\vee
\; = \; \gerI\big(\gerk^\perp\big) \, ,  &  \phantom{\big|}
\quad  \calC(K)_\h^{\,\triangledown} \mod \h \, {F_\h[[G]]}^\vee
\; = \; \gerC\big(\gerk^\perp\big) \, ,  \\
   \gerI(\gerk)_\h^{\;!} \mod \h \, {U_\h(\gerg)}'
\; = \; \calI\big(K^\perp\big) \, ,  &  \phantom{\Big|}
\quad  \gerC(\gerk)_\h^{\,\Lsh} \mod \, \h \, {U_\h(\gerg)}'
\; = \; \calC\big(K^\perp\big) \, .
      \end{matrix}  $$
 \vskip-7pt
\noindent
 That is,  $ \Big( \calI(K)_\h^{\,\curlyvee} , \, \calC(K)_\h^{\,
\triangledown} , \, \gerI(\gerk)_\h^{\;!} \, , \, \gerC(\gerk)_\h^{\,
\Lsh} \Big) $  is a 
         quantization of the quadruple\break
$ \phantom{\bigg|} \! \Big( \gerI\big(\gerk^\perp\big) ,
\, \gerC \big( \gerk^\perp\big) , \, \calI\big(K^\perp\big) ,
\calC \big(K^\perp\big) \!\Big) $  w.r.t.~the quantization
     $ \Big( {F_\h[[G]]}^\vee , {U_\h(\gerg)}' \Big) $\break
of  $ \, \Big( U(\gerg^*) \, , \, F[[G^*]] \Big) \, $,  \,
which again satisfies relations like in  $ (3.2) \, $.}
\end{theorem}   

\noindent
 {\it Sketch of the proof.}  Let us draw a quick sketch of the
proof of Theorem 4.1, as it is given in the extended version
[3].  The main idea is to reduce everything to the study of 
$ \calI(K)_\h^{\,\curlyvee} $  and  $ \calC(K)_\h^{\,\triangledown} $, 
and to get the rest via an indirect approach.  Indeed, this amounts to
show that the quadruple  $ \, \Big( \calI(K)_\h^{\,\curlyvee} \, ,
\, \calC(K)_\h^{\,\triangledown} \, , \, \gerI(\gerk)_\h^{\;!} \; ,
\, \gerC(\gerk)_\h^{\,\Lsh} \Big) \, $  satisfies relations similar
to (3.2), namely
 \vskip-15pt
  $$  \hskip-2pt   \hbox{ $ \begin{matrix}
   {\calI_\h}^{\!\curlyvee} = {\big( {\gerC_\h}^{\!\!\Lsh\,} \big)}^\perp
\, ,  \quad  {\gerC_\h}^{\!\!\Lsh} = {\big( {\calI_\h}^{\!\curlyvee}
\big)}^\perp \, ,  \quad  &  {\gerI_\h}^{\! !} = {\big( {\calC_\h}^{\!\!
\triangledown} \big)}^\perp \, ,  \quad  {\calC_\h}^{\!\!\triangledown}
= {\big( {\gerI_\h}^{\! ! \,} \big)}^\perp  \\
   {\calI_\h}^{\!\curlyvee} \! = {F_\h}^{\!\!\vee} {\big( {\calC_\h}^{\!
\!\triangledown} \big)}^{\!+} ,  \;  {\calC_\h}^{\!\!\triangledown} \! =
{\big( {F_\h}^{\!\vee} \big)}^{\!\text{\it co}{\calI_\h}^{\!\curlyvee}} ,
\hskip3pt  &  {\gerI_\h}^{\! !} \! = {U_\h}^{\!\prime} \, {\big(
{\gerC_\h}^{\!\!\Lsh} \,\big)}^{\!+} ,  \;  {\gerC_\h}^{\!\!\Lsh}
\! = {\big( {U_\h}^{\!\prime} \,\big)}^{\!\text{\it co}{\gerI_\h}^{\! !}}
   \end{matrix} $ }   \hskip3pt (4.1)  $$
 The relations in the top line of (4.1) follow from the similar relations
for the elements of the initial quadruple   --- the top line relations in
(3.2) ---   passing through the duality-preserving property of Drinfeld's
functors given by Proposition 3.1.  In fact, these orthogonality relations are proved much like Proposition 3.1 itself.  Similarly, the bottom line relations in (4.1) follow from the bottom line relations in (3.2) and the very definitions.
                                              \par \noindent
 In force of (4.1), it is enough to prove that just  {\sl one\/}  of the
four objects involved   --- those on left-hand side of the identities in
bottom line of (4.1) ---   has the required properties, in particular, it
is a quantum subgroup for  $ K^\perp $  or a quantum Poisson quotient for
$ \, G^* \big/ K^\perp \, $.  In fact, the similar results for the other
three objects will then follow as a consequence, due to (4.1).

\vskip3pt

\noindent
 Next step is to prove that our quantum sub-objects have the ``right''
Hopf algebraic properties, i.e.~those occurring in first column of
$ (3.1)' $,  namely
  $$  \hskip-3pt   \begin{matrix}
   \phantom{\Big|}  \text{\it (a)} \;\;\; {\calI_\h}^{\!\curlyvee}
\ideal_\ell {F_\h}^{\!\vee}  &  {} \hskip13pt  \text{(\it b)} \;\;\;
{\calC_\h}^{\!\!\triangledown} \leq^1 \! {F_\h}^{\!\vee}
&  {} \hskip13pt  \text{\it (c)} \;\;\; {\gerI_\h}^{\! !}
\ideal_\ell {U_\h}'  &  {} \hskip13pt  \text{\it (d)}
\;\;\; {\gerC_\h}^{\!\!\Lsh} \leq^1 \! {U_\h}'  \\
   \phantom{\Big|}  \text{\it (e)} \;\;\; {\calI_\h}^{\!\curlyvee}
\,\coideal\; {F_\h}^{\!\vee}  &  {} \hskip13pt  \text{\it (f)} \;\;\;
{\calC_\h}^{\!\!\triangledown} \,\coideal_\ell \, {F_\h}^{\!\vee}
&  {} \hskip13pt  \text{\it (g)} \;\;\; {\gerI_\h}^{\! !} \,\coideal\;
{U_\h}'  &  {} \hskip13pt \text{\it (h)} \;\;\; {\gerC_\h}^{\!\!\Lsh}
\,\coideal_\ell\, {U_\h}'
      \end{matrix}  $$
This is an easy task, with some shortcuts available thanks to (4.1)
again.

\vskip3pt

\noindent
 As a third step, we must prove that our quantum sub-objects have
``the good property'' with respect to specialization   --- see
$ (3.1)' $  ---   namely
  $$  \begin{matrix}
   \hskip5pt  \text{\it (a)} \hskip13pt  {\calI_\h}^{\!\curlyvee}
{\textstyle \bigcap} \; \h \, {F_\h}^{\!\vee} \, = \; \h \,
{\calI_\h}^{\!\curlyvee}  &  \qquad  \hskip11pt
   \text{\it (b)} \hskip13pt  {\calC_\h}^{\!\!\triangledown} \,
{\textstyle \bigcap} \; \h \, {F_\h}^{\!\vee} \, = \; \h \,
{\calC_\h}^{\!\!\triangledown}  \\
   \text{\it (c)} \hskip13pt  {\gerI_\h}^{\! !} \, {\textstyle \bigcap}
\; \h \, {U_\h}' \, = \; \h \, {\gerI_\h}^{\! !}  &  \qquad
   \text{\it (d)} \hskip13pt  {\gerC_\h}^{\!\!\Lsh} \,
{\textstyle \bigcap} \; \h \, {U_\h}' \, = \; \h \,
{\gerC_\h}^{\!\!\Lsh}
      \end{matrix}  $$
 Proving this requires a different analysis according to whether one
deals with cases  {\it (a)\/}  and  {\it (b)\/}  or cases  {\it (c)\/}
and  {\it (d)}.  Indeed, the latter identities are proved easily,
just looking at definitions of  $ {\calI_\h}^{\!\curlyvee} $  and
$ {\calC_\h}^{\!\!\triangledown} $  and reminding that  $ \calI_\h $
and  $ \calC_\h $  do satisfy the ``good property'' for specialization,
by assumption.  Instead, cases  {\it (a)\/}  and  {\it (b)\/}  require
a careful description of  $ {\calI_\h}^{\!\curlyvee} $  and
$ {\calC_\h}^{\!\!\triangledown} \, $.
                                              \par \noindent
 Let  $ \, I := I_{\scriptscriptstyle F_\h} \, $  and  $ \, J :=
J_{\scriptscriptstyle F_\h} \, $  be as in Subsection 3.3, and  $ \, J^\vee
:= \h^{-1} J \subset {F_\h}^{\!\vee} \, $.  Then  $ \, J \! \mod \h
\, F_\h = J_{\scriptscriptstyle G} := \text{\it Ker}\, \big(\,\epsilon
\, \colon \, F[[G]] \longrightarrow \Bbbk \,\big) \, $,  \, and  $ \,
J_{\scriptscriptstyle G} \big/ {J_{\scriptscriptstyle G}}^{\!2} =
\gerg^* \, $.  Let  $ \, \{y_1,\dots, y_n\} \, $,  with  $ \, n :=
\dim(G) \, $,  be a  $ \Bbbk $--basis of  $ \, J_{\scriptscriptstyle G}
\big/ {J_{\scriptscriptstyle G}}^{\!2} \, $,  and pull it back to a
subset  $ \, \{j_1,\dots,j_n\} \, $  of  $ J \, $.  Then  $ \, \big\{
\h^{-|\underline{e}|} j^{\,\underline{e}} \mod \h \, {F_\h}^{\!\vee}
\;\big|\; \underline{e} \in \N^{\,n} \,\big\} \, $  (with  $ \, j^{\,
\underline{e}} := \prod_{s=1}^n j_s^{\,\underline{e}(i)} \, $),  is a
$ \Bbbk $--basis  of  $ {F_0}^{\!\vee} $  and, if  $ \, j_s^{\,\vee}
:= \h^{-1} j_s \, $  for all  $ s $,  the set  $ \, \big\{ j_1^{\,\vee},
\dots, j_n^{\,\vee} \big\} \, $  is a  $ \Bbbk $--basis  of  $ \; \gert
:= J^\vee \mod \h \, {F_\h}^{\!\vee} \, $.  Moreover, since  $ \, j_\mu
\, j_\nu - j_\nu \, j_\mu \in \h \, J \, $  (for  $ \, \mu, \nu \in \{
1, \dots, n\} \, $)  we have  $ \; j_\mu \, j_\nu - j_\nu \, j_\mu = \h
\sum_{s=1}^n c_s \, j_s + \h^2 \gamma_1 + \h \, \gamma_2 \; $  for some
$ \, c_s \in \kh \, $,  $ \, \gamma_1 \in J \, $  and  $ \, \gamma_2 \in
J^2 $,  whence  $ \; \big[ j_\mu^\vee, j_\nu^\vee \,\big] := j_\mu^\vee
\, j_\nu^\vee - j_\nu^\vee \, j_\mu^\vee \equiv {\textstyle \sum_{s=1}^n}
\, c_s \, j_s^\vee \; \mod \, \h \, {F_\h}^{\!\vee} \, $,  thus  $ \;
\gert := J^\vee \! \mod \h \, {F_\h}^{\!\vee} \, $  is a Lie subalgebra
of  $ {F_0}^{\!\vee} \, $:  indeed,  $ \, {F_0}^{\!\vee} = U(\gert) \, $
as Hopf algebras.  Even more, this description also shows that the linear
map  $ \, \gert \longrightarrow \gerg^* \, $  given by  $ \, y_s \mapsto
j_s^\vee \; \big(\!\!\! \mod \h \, {F_\h}^{\!\vee} \big) \, $,  $ \, s
= 1, \dots, n \, $,  is a Lie bialgebra isomorphism.
                                              \par \noindent
 Let us now fix the set  $ \, \{y_1, \dots, y_n\} \, $  as follows. 
If  $ \, k := \dim(K) \, $,  we can choose a system of parameters for 
$ G \, $,  say  $ \, \big\{ j_1, \dots, j_k, j_{k+1}, \dots, j_n \big\}
\, $  such that  $ \, \calC(K) := {F[[G]]}^K \! = \Bbbk [[ j_{k+1}, \dots, j_n ]] \, $,  the topological  {\sl subalgebra\/}  of  $ F[[G]] $  generated by  $ \big\{ j_{k+1}, \dots, j_n \big\} \, $,  and  $ \calI(K) = \! \big( j_{k+1}, \dots, j_n \big) \, $,  the  {\sl ideal\/}  of  $ F[[G]] $  topologically generated by  $ \big\{ j_{k+1}, \dots, j_n \big\} \, $.  Set also  $ \, y_s := j_s \! \mod {{J_{\scriptscriptstyle G}}^{\!2}} \; (s = 1, \dots n) \, $.  Then  $ \, \hbox{\it Span}\,\big( \{ y_{k+1}, \dots, y_n\} \big) = \gerk^\perp \, $.
                                              \par \noindent
 Basing on this, one finds that  $ \, \calC(K)_\h^{\,\triangledown} \, $ 
is just the topological subalgebra of  $ \, \uhg = \Bbbk \big[ j_1^\vee,
\dots, j_n^\vee \big][[\h]] \, $  generated by  $ \, \big\{ j_{k+1}^\vee,
\dots, j_n^\vee \big\} \, $,  that is  $ \, \Bbbk \big[ j_{k+1}^\vee,
\dots, j_n^\vee \big][[\h]] \, $.  Similarly,  $ \, {\calI_\h}^{\!\curlyvee}
\, $  is the left ideal of  $ \uhg $  generated by  $ \, \big\{ j_{k+1}^\vee,
\dots, j_n^\vee \big\} \, $,  that is the set of all  series (in  $ \h
\, $)  in  $ \, \Bbbk \big[ j_1^\vee, \dots, j_n^\vee \big] [[\h]] \, $ 
whose coefficients belong to the ideal of  $ \Bbbk \big[ j_1^\vee, \dots,
j_n^\vee \big] $  generated by  $ \, \big\{ j_{k+1}^\vee, \dots, j_n^\vee
\big\} \, $.  Then  $ (3.1)' $  implies that  $ \, \calC(K)_\h^{\,
\triangledown} \cap \h \, {F_\h[G]}^\vee \subseteq \, \h \,
\calC(K)_\h^{\,\triangledown} \, $,  while the converse is obvious. 
This proves  {\it (b)},  and case  {\it (a)\/}  is similar.

\vskip3pt

\noindent
 At this point, one has proved that each element of the quadruple  $ \,
\Big(\, \calI(K)_\h^{\,\curlyvee} \, , \, \calC(K)_\h^{\,\triangledown}
\, , \, \gerI(\gerk)_\h^{\;!} \; , \, \gerC(\gerk)_\h^{\,\Lsh} \,\Big)
\, $  is a quantum subgroup   --- the second and third element ---
or a quantum Poisson quotient   --- the first and fourth element ---
of the dual Poisson group  $ G^* \, $,  with respect to the fixed
quantization  $ \; \Big(\, {F_\h[[G]]}^\vee , \, {U_\h(\gerg)}' \,\Big)
\; $  of  $ \; \Big(\, U(\gerg^*) \, , \, F[[G^*]] \,\Big) \, $.  In
addition, relations (4.1) induce similar relations when specializing
at  $ \, \h = 0 \, $,  so these quantum subgroups and quantum Poisson
quotients all provide a quantization of  {\sl one\/}  single pair
$ \big(\, T \, , \, G^* \!\big/ T \,\big) $,  for some coisotropic
(formal) subgroup  $ T $  of  $ G^* \, $.

\vskip3pt

\noindent
 The fourth step concerns last part of claim  {\it (a)\/}:  it
amounts to prove that
 \vskip-16pt
  $$  {\big( \calI_\h^{\,\curlyvee} \big)}^! \, = \; \calI_\h \quad ,
\;\quad  {\big( \calC_\h^{\,\triangledown} \big)}^\Lsh \, = \; \calC_\h
\quad ,  \;\quad  {\big( {\gerI_\h}^{\! ! \,} \big)}^{\!\curlyvee} \, =
\; \gerI_\h \quad ,  \;\quad  {\big( {\gerC_\h}^{\!\!\Lsh\,} \big)}^{\!
\triangledown} \, = \; \gerC_\h  $$
 \vskip-4pt
\noindent
 Now, the very definitions imply at once one-way inclusions
 \vskip-16pt
  $$  {\big( \calI_\h^{\,\curlyvee} \big)}^! \, \supseteq \; \calI_\h
\quad , \;\quad  {\big( \calC_\h^{\,\triangledown} \big)}^\Lsh \, \supseteq
\; \calC_\h \quad ,  \;\quad  {\big( {\gerI_\h}^{\! ! \,} \big)}^{\!
\curlyvee} \, \subseteq \; \gerI_\h \quad ,  \;\quad  {\big( {\gerC_\h}^{\!
\!\Lsh\,} \big)}^{\!\triangledown} \, \subseteq \; \gerC_\h  $$
 \vskip-4pt
\noindent
 For the converse inclusions, in the first or the second case they follow
again from the description of  $ \, \calI_\h^{\,\curlyvee} $  or  $ \,
\calC_\h^{\,\triangledown} $,  respectively.  Then one uses the identities
just proved, along with the orthogonality-preserving properties of
Drinfeld-like functors   --- the top line in (4.1) ---   applied twice,
to obtain the full identities in the other cases too.  Note that (again)
one could simply prove only one of the identities involved, and then
get the others via (4.1).

\vskip3pt

\noindent
 The last step is to show that the coisotropic (formal) subgroup  $ T $ 
of  $ G^* $  found above is  $ K^\perp $,  i.e.~$ \, \text{\it Lie}\,(T)
= \gerk^\perp \, $.  This means that one has to specialize our quantum
subobjects at  $ \, \h = 0 \, $.  By the third step, it is enough to
do it for any one of them.  The best choice is again 
$ \calC_\h^{\,\triangledown} $  or  $ \calI_\h^{\,\curlyvee} $, 
whose (almost) explicit description yields an explicit
description of its 
   \hbox{specialization too.   $ \, \square $}   

\subsection{Remark}     
 We point out that quantum coisotropic subgroups such
as  $ \calI(K)_\h $  and  $ \gerC(\gerk)_\h $  provide quantum Poisson
quotients  $ \, \calI(K)_\h^{\,\curlyvee} = \gerI \big( \gerk^\perp
\big)_\h \, $  and  $ \, \gerC(\gerk)_\h^{\,\Lsh} = \calC \big( K^\perp
\big)_\h \, $  respectively, while the quantum Poisson quotients  $ \,
\calC(K)_\h \, $  and  $ \, \gerI(\gerk)_\h \, $  yield quantum
coisotropic subgroups  $ \, \calC(K)_\h^{\,\triangledown} = \gerC
\big(\gerk^\perp\big)_\h \, $  and  $ \, \gerI(\gerk)_\h^{\;!} =
\calI\big(K^\perp\big)_\h \, $  respectively.  Thus, Drinfeld-like
functors map quantum coisotropic subgroups to quantum Poisson
quotients, and viceversa.

\vskip0,6truecm

\section{Generalizations and effectiveness}

\vskip5pt

\subsection{QDP with half quantizations}     
 In this work we start from
a pair of mutually dual quantum groups,  i.e.~$ \big( \fhg \, , \,
\uhg \big) \, $.  This is used in the proofs to apply orthogonality
arguments.  However, this is only a matter of choice.  Indeed, our
QDP deals with quantum subgroups or quantum homo\-geneous spaces
which are contained either in  $ \fhg $  or in  $ \uhg \, $.  In
fact  {\sl we might prove every step in our discussion using only
the single quantum group which is concerned},  and only one quantum
subgroup (such as  $ \calI_\h $,  or  $ \calC _\h $,  etc.) at the
time, by a direct method with no orthogonality arguments.   

\subsection{QDP with global quantizations}     
 In this paper we use quantum
groups in the sense of Definition 3.1; these are sometimes called 
{\sl local\/}  quantizations.  Instead, one can consider  {\sl
global quantizations\/}:  quantum groups like Jimbo's, Lusztig's,
etc.  The latter differ from the former in that
 \vskip0pt
\noindent
 \quad   {\it --1)}  \, they are standard (rather than topological)
Hopf algebras;
                                           \par \noindent
 \quad  {\it --2)}  \, they may be defined over any ring  $ R \, $,
the r\^{o}le of  $ \h $  being played by a suitable element of that
ring (for example,  $ \, R = \Bbbk \big[ q, q^{-1} \big] \, $  and
$ \, \h = q - 1 \, $).
 \vskip4pt
\noindent
 Now, our analysis may be done also in terms of  {\sl global quantum
groups\/}  and their specializations.  First, one starts with  {\sl
algebraic\/}  (instead of formal) Poisson groups and Poisson homogeneous
spaces.  Then one defines Drinfeld-like functors in a similar manner;
the key fact is that the QDP for quantum groups has a  {\sl global 
version [12]},  and the recipes of Section 4 to define Drinfeld-like
functors still do make sense.  Moreover, one can extend our QDP for
coisotropic subgroups (and Poisson quotients) to  {\sl all\/}  closed
sub\-groups and homogeneous spaces: all this will be treated separately.  

\subsection{$ * $--structures  and QDP in the real case}     
 If one looks for quantizations of  {\sl real\/}  subgroups and
homogeneous spaces, then one must consider  $ * $--structures 
on the quantum group Hopf algebras.  Then one can perform all our construction in this setting, and state and prove a version of the QDP for  {\sl real\/}  quantum subgroups and quantum homogeneous spaces too, both in the formal and in the global setting.   

\subsection{QDP for pointed Poisson varieties}     
 Let {\sl pointed Poisson variety (=p.P.v.)\/}  be any pair 
$ (M,\bar{m}) $  where  $ M $  is a Poisson variety and  $ \,
\bar{m} \in M \, $  is such that  $ \{\bar{m}\} $  is a symplectic
leaf of  $ M \, $.  A  {\sl morphism of p.P.v.'s  $ (M, \bar{m}) $, 
$ (N,\bar{n}) $\/}  is any Poisson map   $ \, \varphi : M \longrightarrow
N \, $  such that  $ \, \varphi(\bar{m}) = \bar{n} \, $.  This defines a
subcategory of the category of all Poisson varieties, whose morphisms are
those which map distinguished points into distinguished points.  In terms
of affine algebraic geometry, a  p.P.v.~$ (M,\bar{m}) $  is given by the
pair  $ \big( F[M], \mathfrak{m}_{\bar{m}}\big) $  where  $ F[M] $  is
the function algebra of  $ M $  and $ \mathfrak{m}_{\bar{m}} $  is the
defining ideal of  $ \, \bar{m} \in M \, $  in  $ F[M] \, $.
                                               \par \noindent
 By assumptions, the Poisson bracket of  $ F[M] $  restricts to a Lie bracket onto  $ \mathfrak{m}_{\bar{m}} \, $:  then
[16]  $ \, \mathcal{L}_M := \mathfrak{m}_{\bar{m}} \big/
\mathfrak{m}_{\bar{m}}^{\;2} \, $  (the cotangent space to 
$ M $  at  $ \bar{m} $)  inherits a Lie algebra structure too,
the so-called ``linear approximation of  $ M $  at  $ \bar{m} \, $''.
                                         \par \noindent
 Poisson quotients are natural examples of p.P.v.'s.
Other examples are Poisson monoids (\,=\,unital Poisson semigroups),
each one being pointed by its unit element.  If  $ \, (M,\bar{m}) =
(\varLambda,e) \, $  is a Poisson monoid, then  $ F[\varLambda] $
is a bialgebra (and conversely), and  $ \mathcal{L}_\varLambda
\, $  has a natural structure of  {\sl Lie bialgebra},  hence
$ U(\mathcal{L}_\varLambda) $  is a co-Poisson Hopf algebra.  The
Lie cobracket is induced by the coproduct of  $ F[\varLambda] $,
hence (dually) by the multiplication in  $ \varLambda \, $.  In
particular, when the monoid  $ \varLambda $  is a Poisson group
$ G $  we have  $ \, \varLambda = \gerg^* \, $.
                                    \par \noindent
 We call   {\sl quantization of a p.P.v.}~$ (M,\bar{m}) $  any
unital algebra  $ A $  in  $ \mathcal{T}_{\otimeshat} \, $  or 
$ \mathcal{P}_{\otimestilde} $  (see Subsection 3.1) along
with a morphism of topological unital algebras  $ \,
{\underline{\epsilon}}_{{}_A} \colon A \longrightarrow \kh
\, $,  such that  $ \, A{\big|}_{\h=0} \! \cong F[M] \, $  as
Poisson  $ \Bbbk $--algebras  and  $ \, \pi_{{}_{\Bbbk[[\h]]}} \!
\circ {\underline{\epsilon}}_{{}_A} \! = {\underline{\epsilon}}_{{}_M}
\! \circ \pi_{{}_A} \, $,  with  $ \, \pi_{{}_A} \colon A
\! \relbar\joinrel\twoheadrightarrow \! A{\big|}_{\h=0}
\! \cong F[M] \, $,  $ \, \pi_{{}_{\Bbbk[[\h]]}} \colon
\kh \! \relbar\joinrel\twoheadrightarrow \kh \, $  the
specialisation maps  ($ \, \h \mapsto 0 \, $);  in
this case we write  $ \, A = F_\h[M] \, $.  For any
such object we set  $ \, J_{\h,M} := \text{\sl Ker}\,\big(
{\underline{\epsilon}}_{{}_{F_\h[M]}} \big) \, $  and  $ \,
I_{\h,M} := J_{\h,M} + \h \, F_\h[M] \, $.  A  {\sl morphism
of quantizations of p.P.v.'s\/}  is any morphism  $ \, \phi \,
\colon F_\h[M] \! \longrightarrow \! F_\h[N] \, $  in  $ \mathcal{A}^+ $  such that  $ \, {\underline{\epsilon}}_{{}_{F_\h[N]}} \! \circ \phi = {\underline{\epsilon}}_{{}_{F_\h[M]}} \, $.  Quantizations of
p.P.v.'s and their morphisms form a subcategory of 
$ \mathcal{T}_{\otimeshat} \, $,  or  $ \mathcal{P}_{\otimestilde}
\, $,  respectively.  Also, we might repeat the same construction
in the setting of formal geometry, or we might use global
quantizations, as in Subsection 5.2.
                                    \par \noindent
 Now define  $ {F_\h[M]}^\vee $  like in Subsection 3.3, replacing  $ H $ 
with  $ \, A = F_\h[M] \, $  and  $ J_{\scriptscriptstyle H} $  with 
$ J_{\scriptscriptstyle \h, M} \, $.  Then the same analysis made to
prove the parts of Theorem 4.1 concerning  $ {F_\h[[G]]}^\vee $
proves also the following:

\begin{theorem}     
 {\it Let  $ \, F_\h[M] \in \mathcal{A}^+ \, $  be
a quantization of a pointed Poisson manifold  $ (M,\bar{m}) $  as above.
Then  $ {F_\h[M]}^\vee $  is a quantization of
$ \, U(\mathcal{L}_M) $,  i.e.
 \vskip-14pt
  $$  {F_\h[M]}^\vee{\Big|}_{\h=0} \, := \; {F_\h[M]}^\vee
\! \Big/ \h \, {F_\h[M]}^\vee \, = \; U(\mathcal{L}_M)  $$
 \vskip-3pt
\noindent
 If in addition  $ M $  is a Poisson monoid and  $ F_\h[M] $  is
a quantization of  $ F[M] $  {\sl as a bialgebra},  then the last
identification above is one of Hopf algebras.
                                           \par \noindent
 Moreover, the construction  $ \; F_\h[M] \mapsto {F_\h[M]}^\vee \; $ 
is functorial.}
\end{theorem}

\subsection{Computations}     
 In the extended version [3] of this work a nontrivial example
is treated in detail.  Here instead, let us consider the 1--parameter
family of  {\it quantum spheres\/}  $ \mathbb{S}^2_{q,c} $  described by
Dijkhuizen and Koornwinder [4].  These are quantum homogeneous
spaces for the standard  $ SU_q(2) $,  described by objects of type 
{\it (c)}.  Such objects, corresponding to $1$--dimensional subgroups
conjugated to the diagonally embedded  $ \mathbb{S}^1 $,  are right
ideals and two sided coideals in  $ U_q\big(\mathfrak{su}(2)\big) $ 
generated by a single twisted primitive element  $ X_\rho \, $,  i.e. 
$ \, X_\rho = \big( q^2 + 1 \big)^{\! -1/2} L^{-1} E \, + \, \big( q^2
+ 1 \big)^{\! -1/2} F \, L \, - \, \rho \, \frac{\;{( q + \,  q^{-1}
)}^{1/2}}{q \, - \, q^{-1} } \big( L - L^{-1} \big) \, $,  \, using
notations of the work of  Gavarini [12].  Applying the suitable
duality functor, these elements correspond to generators of dual
right ideals and two sided coideals in  $ F_q\big[SU(2)^*\big]
\, $,  say  $ \, \xi_\rho = x + y - \rho \, \big( z - z^{-1}
\big) \, $  (in Gavarini's [12] notation again).  Modding
out these elements provides the corresponding quantum coisotropic
subgroups.  Here much of the Poisson aspects of the theory trivializes. 
It would be interesting to carry out similar computations for the
1--parameter families of quantum projective spaces described by
Dijkhuizen and Noumi [5].  
  
\vskip0,6truecm

\end{document}